\documentclass{article}
\usepackage{amsthm}
\usepackage{amsmath}
\usepackage{graphicx} 
\usepackage{color}

\usepackage[letterpaper,top=2cm,bottom=2cm,left=3cm,right=3cm,marginparwidth=1.75cm]{geometry}

\newtheorem{thm}{Theorem}

\newtheorem{conj}[thm]{Conjecture}

\newtheorem{lemma}[thm]{Lemma}

\bfseries\normalfont

\title{A Sublinear Minimum-Degree Condition for $2$-Connected Subgraphs of All Orders}
	\author{\bf Kenta Ozeki\footnote{
    Yokohama National University, Japan. ozeki-kenta-xr@ynu.ac.jp.} \qquad {\bf Takahiro Ueoro\footnote{
    Yokohama National University, Japan. ueoro-takahiro-wc@ynu.jp.}}}

\begin{document}
\maketitle

\begin{abstract}
Motivated by an analogue of pancyclicity, we study minimum-degree conditions ensuring that a $2$-connected graph $G$ of order $n$ contains a $2$-connected subgraph of every order $\ell\in\{4,5,\ldots,n\}$. Yin and Wu [A minimum degree condition for a 2-connected graph containing all possible orders of 2-connected subgraphs, Discrete Appl. Math. 387 (2026), 129–136]
initiated the study of this problem and showed that the condition $\delta(G)\ge \lceil n/3\rceil+1$ is sufficient. Kashima conjectured that the condition $\delta(G)\ge \sqrt{3n}$ is sufficient. In this paper, we prove that every $2$-connected graph $G$ of order $n$ with $\delta(G)\ge 2n^{2/3}+6n^{1/3}+2$ contains a $2$-connected subgraph of every order from $4$ to $n$. In particular, this gives the first sufficient minimum-degree condition of sublinear order in $n$.
\end{abstract}

\section{Introduction}

Extending the concept to pancyclicity, some researchers have recently considered the existence of a $2$-connected subgraph of order $\ell$ for every $\ell \in \{4,5,\dots,n\}$ in a $2$-connected graph $G$ of order $n$.
Yin and Wu \cite{YW} showed that this property holds if $\delta(G) \geq \lceil n/3 \rceil + 1$. This was subsequently improved by Liu and Ning \cite{LN}, who showed that it holds under the weaker condition $\delta(G) \geq \lceil n/4 \rceil + 2$. The current best known result is the following theorem due to Yang. This result was conjectured in \cite{LN}.

\begin{thm}[Yang, \cite{HY}]
For every fixed integer $q \ge 3$,
there exists an integer $n_0(q)$ such that
every $2$-connected graph $G$ of order $n \ge n_0(q)$ with
$$\delta(G) \ge \frac{n}{q}$$
contains a $2$-connected subgraph of order $\ell$ for every $\ell \in \{4, 5, \dots , n\}$.
\end{thm}

On the other hand, Yin and Wu \cite{YW} conjectured that the same conclusion holds even if $\delta(G) \ge \sqrt{n}$. This conjecture was subsequently disproved by Liu and Ning \cite{LN}.
Kashima \cite{Kashima} pointed out the following construction of a $2$-connected graph $G$ of order $n$ with $\delta(G) = \sqrt{3n+1}-2$ that does not contain a $2$-connected subgraph of every order. Let $t \ge 2$, and let $H$ be the graph obtained from $K_{3t-2}$ by deleting one edge. Take $t$ copies $H_1,\dots,H_t$ of $H$, and let $x_i$ and $y_i$ be the vertices of $H_i$ corresponding to the endpoints of the deleted edge. Let $G$ be the graph obtained by adding the $t$ edges $x_i y_{i+1}$ for $1 \le i \le t$, where $y_{t+1}$ is understood to be $y_1$.
Note that
$$
n=|V(G)|=t(3t-2) \quad\text{and}\quad \delta(G)=3t-3=\sqrt{3n+1}-2.
$$
However, $G$ contains no $2$-connected subgraph of order $3t-1$. Indeed, any $2$-connected subgraph meeting more than one $H_i$ must meet every $H_i$ in at least three vertices, and hence has order at least $3t$. 

Motivated by this construction, Kashima conjectured the following.

\begin{conj}[Kashima \cite{Kashima}]
Every $2$-connected graph $G$ of order $n$ with
$\delta(G) \ge \sqrt{3n}$
contains a $2$-connected subgraph of order $\ell$ for every
$\ell \in \{4,5,\dots,n\}$.
\end{conj}

In this paper, we give a partial solution to his conjecture as follows.
In particular, this is the first result to guarantee the desired property under a sublinear minimum-degree condition.

\begin{thm}
\label{main}
Let $G$ be a $2$-connected graph of order $n$.
If 
$$\delta(G) \geq 2n^{2/3}+ 6n^{1/3} + 2,$$
then $G$ contains a $2$-connected subgraph of order $\ell$
for every $\ell \in \{4,5,\dots,n\}$.
\end{thm}

We emphasize that our main focus is on the order $n^{2/3}$
of the minimum-degree condition.
We do not attempt to optimize the coefficient or the lower-order terms;
they are chosen for the simplicity of the proof.
In the next section, we give some preliminary with some lemmas, and then prove Theorem \ref{main} in Section \ref{proof_sec}.

\section{Preliminary}

In our proof, subgraphs isomorphic to $K_{2,t}$ play a central role.
To find such a subgraph, we use the following theorem,
which can be directly obtained from the upper bound for the Tur\'{a}n number of $K_{2,t}$;
see \cite{Furedi, KST}.
The bound in our theorem below is somewhat weaker than the bound obtained from the Tur\'{a}n number.
However, we use this weaker bound for the simplicity of the arguments that follow, since it has the same order.
For completeness, we include a proof.
We call a path of length $\ell$ an \emph{$\ell$-path}.

\begin{thm}
\label{K2t_thm}
Let $t$ be an integer with $t\ge 2$, and let $H$ be a graph of order $n$.
If
\[
\delta(H)\ge \sqrt{tn}+1,
\]
then $H$ contains a subgraph isomorphic to $K_{2,t}$.
\end{thm}

\begin{proof}
Suppose, to the contrary, that $H$ contains no subgraph isomorphic to
$K_{2,t}$.
We count the number of $2$-paths in $H$.
For each vertex $x\in V(H)$,
the number of $2$-paths with $x$ as the internal vertex is
$\binom{\deg_H(x)}{2}$.
On the other hand, for each pair of vertices $u,v\in V(H)$,
there are at most $t-1$ $2$-paths with $u$ and $v$ as the end vertices,
since otherwise $u$ and $v$ together with $t$ common neighbors
would form a copy of $K_{2,t}$.
Thus, 
\begin{align*}
0 
&\le (t-1)\binom{n}{2} - \sum_{x\in V(H)} \binom{\deg_H(x)}{2} \\
&\le (t-1)\binom{n}{2} - n \binom{\delta(H)}{2} \\
&< (t-1) \cdot \frac{n^2}{2} - n \cdot \frac{tn}{2} \\
&= - \frac{n^2}{2} < 0, 
\end{align*}
which is a contradiction. 
\end{proof}

Faudree, Ordman, Schelp, Jacobson, and Tuza \cite{FOSJT} showed that,
in a $k$-connected graph $H$ with two specified vertices $u$ and $v$,
there exist $m$ edge-disjoint paths connecting $u$ and $v$, where $m\le k$,
such that the length of each path is bounded in terms of
$|V(H)|$ and $\delta(H)$; see also \cite{GY}.
Using a similar idea, we prove the following lemma,
which will be used in our proof.

\begin{lemma}
\label{pathsum_lemma}
Let $G$ be a $2$-connected graph,
and let $X,Y\subseteq V(G)$ be disjoint sets with
$|X|,|Y|\ge2$.
Then there exist two vertex-disjoint paths $P_1$ and $P_2$
such that, for each $i\in\{1,2\}$, $P_i$ connects a vertex of $X$
to a vertex of $Y$ and contains no other vertex of $X\cup Y$, and
\[
|V(P_i)|\le
\frac{3|V(G)|-6}{\delta(G)}.
\]
\end{lemma}

\begin{proof}
If $\delta(G) = 2$, then the conclusion is immediate.
Thus, we may assume that $\delta(G) \geq 3$.

By Menger's theorem, there exist two vertex-disjoint paths $P_1$ and $P_2$
such that, for each $i\in\{1,2\}$, $P_i$ connects a vertex of $X$
to a vertex of $Y$ and contains no other vertex of $X\cup Y$.
We choose such paths so that $|V(P_1)| + |V(P_2)|$ is as small as possible.
We now only show that $|V(P_1)| \leq \frac{3|V(G)| - 6}{\delta(G)}$
since by symmetry, we can obtained the same upper bound on $|V(P_2)|$.

For $i\in\{1,2\}$, let
\[
P_i=x_i^1,x_i^2,\dots,x_i^{p_i},
\]
where $x_i^1\in X$, $p_i = |V(P_i)|$ and $x_i^{p_i}\in Y$.
Let 
$$S_1 = \{ x_1^j : 1 \le j \le p_1 \text{ and } j \equiv 1 \pmod{3}\}.$$
By definition, $|S_1|=\left\lceil p_1/3 \right\rceil$. 
In particular, $S_1 \neq \emptyset$. 
By the minimality of $|V(P_1)|+|V(P_2)|$, $P_1$ has no chord, and hence $S_1$ is independent. 
Furthermore, we see that 
$$e(S_1, V(P_1)) \le 2|S_1| -1.$$

We first claim that, for any $x_1^j,x_1^{j'}\in S_1$,
$$
N_G(x_1^j)\cap N_G(x_1^{j'})
\setminus \bigl(V(P_1)\cup V(P_2)\bigr)=\emptyset.
$$
Indeed, otherwise, say $1\le j<j'\le p_1$, and let
$v\in N_G(x_1^j)\cap N_G(x_1^{j'})
\setminus \bigl(V(P_1)\cup V(P_2)\bigr)$.
Then
$$
\begin{cases}
    x_1^1,\dots,x_1^j,v,x_1^{j'},\dots,x_1^{p_1} & \text{if $v \notin X \cup Y$,} \\
    v,x_1^{j'},\dots,x_1^{p_1} & \text{if $v \in X$, and} \\
    x_1^1,\dots,x_1^j,v & \text{if $v \in Y$} 
\end{cases}
$$
is a path in $G - V(P_2)$ connecting a vertex in $X$ and a vertex in $Y$ but shorter than $P_1$, a contradiction.
Thus, the claim holds. It implies
$$e\left(S_1, V(G) \setminus V(P_1) \cup V(P_2)\right) \leq |V(G)| - p_1 - p_2.$$

We next claim that, for any $x_1^j,x_1^{j'}\in S_1$ with $1 \le j < j' \le p_1$, if $x_1^jx_2^k \in E(G)$, then $x_1^{j'} x_2^{k'} \notin E(G)$ for any $1 \le k' < k$;
since otherwise, say $x_1^{j'} x_2^{k'} \in E(G)$ for some $1 \le k' < k$, then 
$P_1' = x_1^1,\dots,x_1^j,x_2^{k},\dots,x_2^{p_2}$
and
$P_2' = x_2^1,\dots,x_2^{k'},x_1^{j'},\dots,x_1^{p_1}$
are two vertex-disjoint paths with $|V(P_1')| + |V(P_2')| < |V(P_1)| + |V(P_2)|$, a contradiction.
This implies that 
$$e(S_1, V(P_2)) \leq |S_1| + p_2 -1.$$

By putting the above all together, 
we have
\begin{align*}
\delta(G) |S_1| 
&\le \sum_{x \in S_1} \deg_G(x) \\
&\leq 2|S_1|-1 + |V(G)| - p_1 - p_2 + |S_1| + p_2 -1 \\
&= |V(G)| - p_1 + 3|S_1| -2. 
\end{align*}
Since $|S_1| = \left\lceil p_1/3 \right\rceil \ge p_1/3$,
we have 
$$
|V(G)| - 2 \ge (\delta(G)-3)|S_1| + p_1
\ge \frac{1}{3}p_1 \delta(G), 
$$
and hence 
\begin{align*}
|V(P_1)| = p_1 
\le \frac{3|V(G)|-6}{\delta(G)}.
\end{align*}
This completes the proof of Lemma \ref{pathsum_lemma}.
\end{proof}

\section{Proof of Theorem \ref{main}}
\label{proof_sec}

\begin{proof}
Let $G$ be a $2$-connected graph of order $n$ with 
$$\delta(G) \geq 2n^{2/3}+ 6n^{1/3} + 2.$$
Since $$ n-1 \ge \delta(G) \ge 2n^{2/3}+6n^{1/3}+2 \ge 10, $$ 
we have $n\ge 11$. 

We find a $2$-connected subgraph of order $\ell$ by induction on $\ell$.
Case 1 provides the base case, while the main argument is given in Case 2.
\\

\noindent
\textbf{Case 1:} $4 \le \ell \leq \lfloor 4n^{1/3}\rfloor + 2$.

Note that $\delta(G) \ge 2n^{2/3} + 1$.
By Theorem \ref{K2t_thm}, $G$ has a subgraph $Y$ isomorphic to $K_{2,t}$, where $t = \lfloor 4n^{1/3}\rfloor$.
Then by removing $t - \ell + 2$ vertices of degree two from $Y$,
we obtain a $2$-connected subgraph of $G$ of order $\ell$.
\\

\noindent
\textbf{Case 2:} $\lfloor 4n^{1/3} \rfloor +3\le \ell \leq n$.

Let $m = \lfloor 3 n^{1/3}\rfloor + 1$. 
Note that, since $n\ge 8$,
\[
\ell-m
\ge \left\lfloor 4n^{1/3}\right\rfloor + 3 - \left\lfloor 3n^{1/3}\right\rfloor -1
\ge \left\lfloor 3n^{1/3} + 2\right\rfloor - \left\lfloor 3n^{1/3}\right\rfloor + 2
= 4.
\]
In this case, we will show that if $G$ has a $2$-connected subgraph of order $\ell - m$, then so does a $2$-connected subgraph of order $\ell$.
Together with Case 1, this completes the proof.

Let $X$ be a $2$-connected subgraph of $G$ of order $\ell - m$.
Starting with $S=\emptyset$, repeatedly add to $S$ one vertex
$v\in V(G)\setminus\bigl(V(X)\cup S\bigr)$
satisfying $e\bigl(v,V(X)\cup S\bigr)\ge 2m$,
until no such vertex remains.
If $|S| \ge m$, then adding the first $m$ vertices in $S$ to $X$,
we obtain $2$-connected subgraph of $G$ of order $\ell$.
Thus, we may assume $|S| < m$.

Let $H = G - (V(X) \cup S)$.
Note that 
$|V(H)| \le n$.
By the definition of $S$,
we have 
$$\delta(H) \ge \delta(G) - (2m-1) 
\ge (2n^{2/3} + 6n^{1/3} + 2) - 2(3n^{1/3}+1) +1
\ge 2n^{2/3} +1.$$
Thus, by Theorem \ref{K2t_thm}, $H$ has a subgraph $Y$ isomorphic to $K_{2,t}$, where $t = \lfloor 4 n^{1/3}\rfloor$.

By applying Lemma \ref{pathsum_lemma} to $G$,
there exist two vertex-disjoint paths $P_1$ and $P_2$ such that, for each $i\in \{1,2\}$, $P_i$ connects a vertex in $X$ to a vertex of $Y$ and contains no other vertex of $X \cup Y$, and
$$|V(P_i)| \leq \frac{3|V(G)| - 6}{\delta(G)} 
< \frac{3n}{2 n^{2/3}} = \frac{3}{2} n^{1/3}.$$
Let $Z=X\cup P_1\cup P_2\cup Y$.
Since $X$ and $Y$ are $2$-connected and $P_1,P_2$
are vertex-disjoint paths connecting $X$ and $Y$, $Z$ is $2$-connected.
Note that 
\begin{align*}
    |V(Z)| &\ge |V(X)| + |V(Y)| \\
    &= (\ell-m) + 2+t \\
    &= \ell - (\lfloor 3 n^{1/3}\rfloor +1) + 2 + \lfloor 4 n^{1/3}\rfloor \\
    &\ge \ell, \\
\text{and }   |V(Z)| 
&\le |V(X)| + |V(P_1)| + |V(P_2)| + |V(Y)| -4 \\
&\le \ell - (\lfloor 3 n^{1/3}\rfloor +1) + 2 \cdot \frac{3}{2}n^{1/3} + (t+2) -4 \\
&\le  \ell + t -2.
\end{align*}
Thus, $0 \le |V(Z)| - \ell \le t-2$.
Let $Y'$ be the set of vertices of degree two in $Y$
that are not end vertices of $P_i$ for any $i\in\{1,2\}$.
Since $|Y'| \ge t-2$, we can select $|V(Z)|-\ell$ vertices from $Y'$.
Deleting these vertices from $Z$ results in a graph of order $\ell$.
Moreover, the resulting graph is $2$-connected,
since every deleted vertex has degree two in $Z$,
with both neighbors being the two vertices of $Y$ of degree $t$,
which remain in the resulting graph.

This completes the proof of Theorem \ref{main}.
\end{proof}

\section*{Acknowledgements} 
We would like to thank Masaki Kashima for finding an interesting example and for discussions that initiated this study.
Kenta Ozeki was supported by JSPS KAKENHI, Grant Numbers 26K00616 and 26K00618. 
Takahiro Ueoro was supported by YNU-SPRING Program from Yokohama National University, Grant Numbers 74W8110001 and 74W8110002.

\end{document}